\RequirePackage{rotating}
\documentclass[smallextended]{svjour3}
\pdfoutput=1
\usepackage[utf8]{inputenc}
\usepackage{amsmath,amssymb,mathtools}
\allowdisplaybreaks

\usepackage{graphicx}
\usepackage[caption=false, font=footnotesize]{subfig}

\usepackage{booktabs}
\usepackage{natbib} % author-year citations inline with \citet 
\usepackage{rotating}
\usepackage{multirow}
\usepackage{longtable}
\usepackage{array}
\usepackage{hyperref} % Good for drafting and making sure evething references correctly - comment for final edition if unwanted.

\usepackage{dcolumn}
\newcolumntype{T}{>{\textfont0=\the@{.}{.}{1}}c<{\DC@end}}
\newcolumntype{d}{D{.}{.}{1}}
\makeatletter
\def\DC@endright{$\hfil\egroup\@dcolcolor\box\z@\box\tw@\dcolreset}
\def\dcolcolor#1{\gdef\@dcolcolor{\color{#1}}}
\def\dcolreset{\dcolcolor{black}}
\dcolcolor{black}
\makeatother

\newcolumntype{L}[1]{>{\raggedright\let\newline\\\arraybackslash\hspace{0pt}}m{#1}}
\newcolumntype{C}[1]{>{\centering\let\newline\\\arraybackslash\hspace{0pt}}m{#1}}
\newcolumntype{R}[1]{>{\raggedleft\let\newline\\\arraybackslash\hspace{0pt}}m{#1}}

\graphicspath{{figs/}}
\usepackage{xargs}                      % Use more than one optional parameter in a new commands
\usepackage[pdftex,dvipsnames]{xcolor}  % Coloured text etc.
\definecolor{light-gray}{gray}{0.5}
\usepackage[colorinlistoftodos,prependcaption,textsize=tiny]{todonotes}
\newcommandx{\unsure}[2][1=]{\todo[linecolor=red,backgroundcolor=red!25,bordercolor=red,#1]{#2}}
\newcommandx{\change}[2][1=]{\todo[linecolor=blue,backgroundcolor=blue!25,bordercolor=blue,#1]{#2}}
\newcommandx{\info}[2][1=]{\todo[linecolor=OliveGreen,backgroundcolor=OliveGreen!25,bordercolor=OliveGreen,#1]{#2}}
\newcommandx{\improvement}[2][1=]{\todo[linecolor=Yellow,backgroundcolor=Yellow!25,bordercolor=Yellow,#1]{#2}}
\newcommandx{\thiswillnotshow}[2][1=]{\todo[disable,#1]{#2}}

\renewcommand{\L}{\mathcal{L}}
\newcommand{\N}{\mathcal{N}}
\newcommand{\G}{\mathcal{G}}
\newcommand{\F}{\L^{\mathrm{fl}}}
\newcommand{\PA}{\L^{\mathrm{pa}}}
\newcommand{\K}{\mathcal{K}}
\newcommand{\R}{\mathbb{R}}

\newcommand{\reals}{\mathbb{R}}
\newcommand{\cmplx}{\mathbb{C}}
\newcommand{\symm}{\mathbb{S}}
\newcommand{\herm}{\mathbb{H}}
\newcommand{\cone}{\mathcal{K}}
\newcommand{\qcone}{\mathcal{K}_{\mathrm{q}}}
\newcommand{\lb}{\mathrm{l}}
\newcommand{\ub}{\mathrm{u}}

\newcommand{\abs}[1]{\left\vert #1 \right\vert}

\newcommand{\Tr}{\mathbf{tr}}

\newcommand{\eg}{\emph{e.g.}}
\newcommand{\ie}{\emph{i.e.}}

\newcommand{\myTopTrim}{30}
\newcommand{\myBotTrim}{0}
\newcommand{\myLeftTrim}{0}
\newcommand{\myRightTrim}{32}

\begin{document}
%\listoftodos
%\newpage
\title{On the Robustness and Scalability of Semidefinite Relaxation for Optimal Power Flow Problems}
\titlerunning{On the Robustness and Scalability of SDR for OPF Problems}

\author{Anders~Eltved \and
Joachim Dahl \and
Martin~S.~Andersen
%John~Doe,~\IEEEmembership{Fellow,~OSA,}
%and~Jane~Doe,~\IEEEmembership{Life~Fellow,~IEEE}% <-this % stops a space
%\thanks{M. Shell was with the Department of Electrical and Computer Engineering, Georgia Institute of Technology, Atlanta, GA, 30332 USA e-mail: (see http://www.michaelshell.org/contact.html).}% <-this % stops a space
%\thanks{J. Doe and J. Doe are with Anonymous University.}% <-this % stops a space
%\thanks{Manuscript received April 19, 2005; revised August 26, 2015.}
}
\institute{A.~Eltved and M.~S.~Andersen \at Department of Applied Mathematics and Computer Science, Technical University of Denmark, 2800 Kgs.\ Lyngby, Denmark (e-mail: \{aelt,mskan\}@dtu.dk) 
	\and
	J.~Dahl \at MOSEK ApS, Fruebjergvej 3, Symbion Science Park, 2100 Copenhagen, Denmark (e-mail: dahl.joachim@gmail.com)}%

%\markboth{Journal of \LaTeX\ Class Files,~Vol.~14, No.~8, August~2015}{Shell \MakeLowercase{\textit{et al.}}: Bare Demo of IEEEtran.cls for IEEE Journals}

%\listoftodos[Notes]
%\newpage

\maketitle
{
\fontsize{8}{1}
\noindent \textit{This is a pre-print of an article published in Optimization and Engineering. The final authenticated version is available online at:\\ https://doi.org/10.1007/s11081-019-09427-4}
}
\vspace{5mm}
\begin{abstract}
Semidefinite relaxation techniques have shown great promise for nonconvex optimal power flow problems. However, a number of independent numerical experiments have led to concerns about scalability and robustness of existing SDP solvers. To address these concerns, we investigate some numerical aspects of the problem and compare different state-of-the-art solvers. Our results demonstrate that semidefinite relaxations of large problem instances with on the order of 10,000 buses can be solved reliably and to reasonable accuracy within minutes. Furthermore, the semidefinite relaxation of a test case with 25,000 buses can be solved reliably within half an hour; the largest test case with 82,000 buses is solved within eight hours. We also compare the lower bound obtained via semidefinite relaxation to locally optimal solutions obtained with nonlinear optimization methods and calculate the optimality gap. 

\keywords{AC Optimal Power Flow \and Semidefinite Relaxation \and Optimization \and Numerical Analysis}
\end{abstract}

\section{Introduction}
The alternating current optimal power flow (ACOPF) problem is a nonlinear optimization problem that is concerned with finding an optimal operating point for a power system network. Today, more than 60 years after it was first studied by \citet{Carpentier1962}, the problem still receives considerable attention because of the challenging nature of the problem and its important role in power system planning and operation. Many optimization methods have been applied to the ACOPF problem, including general nonlinear optimization techniques, interior-point methods, and meta-heuristic optimization methods \citep{taylor}.

Following the work of \citet{Jabr2006} and \citet{Bai2008}, the use of convex relaxation techniques applied to the ACOPF problem has been explored extensively; see \eg\ \citep{Low:14a,Low:14b} for a recent survey. The interest in these techniques is driven by the fact that the solution to a relaxed problem provides either a globally optimal solution to the original problem or a global lower bound that can be used to assess the quality of locally optimal solutions found by other means. Moreover, a solution to an SDR may also be used to guide a load flow study \citep{Mak:2018} in order to find a feasible operating point. 

Different convex relaxations of the ACOPF problem have been proposed and studied, including a second-order cone relaxation (SOCR) \citep{Jabr2006}, a semidefinite relaxation (SDR) \citep{Bai2008,Lavaei2012}, moment relaxations \citep{MoH:14,JMPG:15}, and more recently, a quadratic convex relaxation (QCR) \citep{coffrin2016,Hijazi2017}. The different relaxations vary in tightness and computational cost. For example, the SDR is generally tighter than the SOCR, but it is generally also more computationally demanding. In an attempt to address the computational cost associated with the SDR, \citet{AHV:14} and \citet{BAL:18} have proposed simpler, weaker SDRs that are cheaper to solve than the standard SDR. The QCR is generally neither weaker nor stronger than the SDR, but it is computationally cheaper and often provides a lower bound of similar quality as that of the SDR.

The high computational cost of solving an SDR of a large ACOPF problem has given rise to concerns about robustness and scalability \citep{Hijazi2016,Hijazi2017,Madani2017}. These concerns are supported by numerical experiments that show that solving the SDR is not only much slower than other approaches, but also more unreliable \citep{coffrin2016}. Our goal with this paper is to address concerns regarding robustness and scalability by demonstrating numerically that an SDR of the ACOPF problem can be solved both reliably and within minutes using commodity hardware, even for large networks with on the order of 10,000 buses. Our contribution is therefore confined to numerical considerations and implementation details (Section \ref{sec:method}) as well as numerical experiments (Section \ref{sec:results}) with the purpose of investigating scalability, accuracy, and robustness for different solvers. What differentiates our implementation from most implementations that have been described and investigated in the literature is the fact that we construct the SDR manually without the use of modeling tools such as YALMIP \citep{Lofberg2004} and CVX \citep{gb08}. Although this manual approach can be both inflexible and cumbersome, it is typically much faster and allows us to control the exact problem formulation, avoiding automatic transformations that may adversely affect the size and conditioning of the SDR problem. We remark that some modeling tools allow some degree of control over the problem formulation (\eg, through options), but it is generally difficult for non-expert users to predict the final problem formulation.

\paragraph*{Notation}\ The set $\qcone^n = \{ (t,x) \in \reals\times \reals^{n-1} \,|\, \|x\|_2\leq t \}$ denotes the second-order cone in $\reals^n$, $\symm^n$ denotes the set of symmetric matrices of order $n$, and $\herm^n$ is the set of Hermitian matrices of order $n$. The sets $\symm_+^n$ and $\herm_+^n$ are the cones of positive semidefinite matrices in $\symm^n$ and $\herm^n$, respectively. Since the symmetric matrices of order $n$ form a vector space of dimension $n(n+1)/2$, the cone $\symm_+^n$ can be reparameterized as $\cone_\mathrm{s}^n = \{ \mathbf{svec}(X) \,|\, X \in \symm_+^n \} \subset \reals^{n(n+1)/2}$ where $\mathbf{svec}(\cdot)$ is an injective function that maps a symmetric matrix of order $n$ to a vector of length $n(n+1)/2$. Similarly, we define $\cone_{\mathrm{h}}^n = \{ \mathbf{hvec}(X)\,|\, X \in \herm_+^n \} \subset \reals^{n^2}$ where $\mathbf{hvec}(\cdot)$ maps a Hermitian matrix of order $n$ to a vector of length $n^2$. The inner product between two matrices $A,B \in \herm^n$ is $\Tr(A^HB)$ where $\Tr(A)$ denotes the trace of a square matrix $A$. Given a complex number $c = a + \jmath b$ where $\jmath = \sqrt{-1}$, $\Re(c)$ denotes the real part $a$, $\Im(c)$ denotes the imaginary part $b$, and $c^*$ denotes the complex conjugate of $c$.

\section{Method} \label{sec:method}
\subsection{The AC Optimal Power Flow Problem}
An AC power system in steady state can be modeled as a directed graph where the set of nodes $\N = \{1,2,\ldots,n\}$ corresponds to a set of $n$ power buses, and the set of edges $ \L \in \N\times \N$ corresponds to transmission lines, \ie, $ (k,l) \in \L $ if there is a line from bus $k$ to bus $l$. The set
$\F \subseteq \L$ consists of all transmission lines with a flow constraint, $\PA \subseteq \L$ consists of all transmission lines with a phase-angle difference constraint, $ \G_k $ denotes a (possibly empty) set of generators associated with bus $k$, and $ \G = \bigcup_{k \in \N} \G_k$ is the set of all generators. The power produced by generator $g \in \G$ is $s_g = p_g + \jmath q_g $, and at each power bus $k \in \N$, we define a complex load (\ie, demand) $ S_k^d = P_k^d + \jmath Q_k^d $, a complex voltage $ v_k $, and a complex current $ i_k $. To simplify notation, we define a vector of voltages $v = (v_1,v_2,\ldots,v_n)$ and a vector of currents $ i = (i_1,i_2,\ldots,i_n) $. With this notation, the ACOPF problem can be expressed as
\begin{subequations} \label{eq:opf}
	\begin{equation}
		\mbox{minimize} \ \ \sum_{g\in\G} f_g(p_g) \
	\end{equation}
	\mbox{subject to} \\
	\begin{align}
							i_k^* v_k & = \sum_{g \in \G_k} s_g - S_k^d, & k&\in \N \label{eq:power_balance} \\
						 	P_g^{\min}	& \le p_g \le P_g^{\max}, & g&\in \G \label{eq:opf:realpower}\\
						 	Q_g^{\min}	& \le q_g \le Q_g^{\max}, & g&\in \G\label{eq:opf:reactivepower} \\
						 	V_k^{\min}	& \le \abs{v_k} \le V_k^{\max}, & k&\in \N \label{eq:opf:voltage}\\
						   \abs{S_{k,l}^{\text{fl}} (v)} & \le S_{k,l}^{\max}, & (k,l)&\in\F	\label{eq:flow_con1}\\
						   \abs{S_{l,k}^{\text{fl}} (v)} & \le S_{l,k}^{\max}, & (k,l)&\in\F \label{eq:flow_con2} \\
						 \phi_{k,l}^{\min} & \le \angle (v_k v_l^*) \le \phi_{k,l}^{\max}, & (k,l) &\in \PA \label{eq:opf:phaseanglediff} \\
						 i &= Yv
						 \label{eq:opf:ohm}
	\end{align}
\end{subequations}
with variables $i \in \cmplx^n$, $v \in \cmplx^n$, and $s \in \cmplx^{|\G|}$, and where $i=Yv$ corresponds to Ohm's law in matrix form, given the network admittance matrix $Y \in \cmplx^{n \times n}$. The cost of generation for generator $g$ is given by $f_g(p_g)$, and we will restrict our attention to convex quadratic generation cost functions, \ie,
\begin{equation}
	 f_g(p_g) = \alpha_g p_g^2 + \beta_g p_g + \gamma_g,
\end{equation}
where the parameters $ \alpha_g \ge 0 $, $ \beta_g $, and $\gamma_g$ are given.
The constraints \eqref{eq:power_balance} are power balance equations, \eqref{eq:opf:realpower} and \eqref{eq:opf:reactivepower} are generation limits, \eqref{eq:opf:voltage} are voltage magnitude limits, \eqref{eq:flow_con1} and \eqref{eq:flow_con2} are transmission line flow constraints, and \eqref{eq:opf:phaseanglediff} are phase-angle difference constraints. The flow from bus $k$ to bus $l$ is given by $ S_{k,l}^{\text{fl}} (v) = v^H T_{k,l} v + \jmath v^H \widetilde{T}_{k,l}v $ (provided that $(k,l) \in \F$ or $(l,k) \in \F$) where $ T_{k,l} \in \herm^n$ and $ \widetilde{T}_{k,l} \in \herm^n$ are given.

\subsection{Semidefinite Relaxation}
Roughly following the steps described in \citep{AHV:14}, we start by reformulating the ACOPF problem \eqref{eq:opf}. Specifically, we perform the following steps:
\begin{enumerate}
	\item Eliminate $i = Yv$ and substitute $P_g^{\min} + p_g^{\lb}$ for $p_g$, $Q_g^{\min} + q_g^{\lb}$ for $q_g$, and $X$ for $vv^H$.
	\item Drop constant terms in the objective:
	\begin{align*}
		f(p_g) &= \alpha_g(P_g^{\min} + p_g^{\lb})^2 + \beta_g (P_g^{\min} + p_g^{\lb}) + \gamma_g\\
		 &= \alpha_g (p_g^{\lb})^2 + \tilde \beta_g p_g^{\lb} + \text{const.}
	\end{align*}
	where $\tilde \beta_g = (\beta_g + 2\alpha_gP_g^{\min})$.
	\item Introduce an auxiliary variable $t_g$ for each $g \in \G^{\text{quad}} = \{g\in\G \,|\, \alpha_g > 0 \}$ and include epigraph constraint
	\[ \alpha_g(p_g^{\lb})^2 \leq t_g \ \Leftrightarrow \begin{bmatrix} 1/2 + t_g\\1/2-t_g \\ \sqrt{2\alpha_g} p_g^{\lb} \end{bmatrix} \in \qcone^3. \]
	\item Introduce slack variables to obtain a standard-form formulation.
\end{enumerate}
These steps yield the equivalent problem
\begin{subequations} \label{eq:opf-reform} \allowdisplaybreaks
	\begin{equation}
		\mbox{minimize} \ \ \sum_{g\in\G} \tilde \beta_g  p_g^{\lb} + \sum_{g\in\G^{\text{quad}}} t_g
	\end{equation}
	\mbox{subject to} \\
	\begin{align}
		\Tr(Y_k X) & = \sum_{g \in \G_k} (P_g^{\min} + p_g^{\lb}) - P_k^d, &  k &\in \N \label{eq:opf-reform-ls} \\
		\Tr(\widetilde Y_k X) & = \sum_{g \in \G_k} (Q_g^{\min} + q_g^{\lb}) - Q_k^d, & k &\in \N \\
	 	p_g^{\lb} + p_g^{\ub} &=  P_g^{\max}-P_g^{\min}, & g &\in \G \\
		q_g^{\lb} + q_g^{\ub} &=  Q_g^{\max}-Q_g^{\min}, & g &\in \G \\
	 	 X_{kk} - \nu_k^{\lb}	&= (V_k^{\min})^2,  & k &\in \N \\
		 X_{kk} + \nu_k^{\ub} &= (V_k^{\max})^2,  & k &\in \N \\
		 z_{k,l} &= \begin{bmatrix} S_{k,l}^{\max} \\ \Tr(T_{k,l} X) \\ \Tr(\widetilde{T}_{k,l} X)\end{bmatrix}, & (k,l)&\in\F	\\
		 z_{l,k} &= \begin{bmatrix}	S_{l,k}^{\max} \\ \Tr(T_{l,k} X) \\ \Tr(\widetilde{T}_{l,k} X)\end{bmatrix}, & (k,l)&\in\F\\
  	w_g &= \begin{bmatrix}1/2 + t_g \\ 1/2- t_g \\ \sqrt{2\alpha_g} p_g^{\lb}\end{bmatrix}, & g &\in \G^{\text{quad}} \\
		\Im(X_{kl}) &=	\tan(\phi_{k,l}^{\min}) \Re(X_{kl}) + y_{k,l}^{\lb}, & (k,l) &\in \PA  \\
	 	\Im(X_{kl}) &=  \tan(\phi_{k,l}^{\max}) \Re(X_{kl}) - y_{k,l}^{\ub}, & (k,l) &\in \PA \label{eq:opf-reform-le} \\
	 p_g^{\lb}, p_g^{\ub} & \ge 0, & g&\in\G \\
	 q_g^{\lb}, q_g^{\ub} & \ge 0, & g&\in\G \\
	\nu_k^{\lb}, \nu_k^{\ub} & \ge 0, & k&\in\N \\
	y_{k,l}^{\lb}, y_{k,l}^{\ub} &\geq 0, & (k,l) &\in \PA\\
	z_{k,l},z_{l,k} &\in \qcone^3, & (k,l) &\in \F \\
  w_g & \in \qcone^3, & g &\in \G^{\text{quad}} \\
	X & = vv^H \label{eq:opf-reform-rank1}
	\end{align}
\end{subequations}
with variables $p^{\lb},p^{\ub},q^{\lb},q^{\ub} \in \reals^{|\G|}$, $t\in \R^{|\G^{\text{quad}}|}$, $\nu^{\lb}, \nu^{\ub} \in \reals^{|\N|}$,
$y^{\lb},y^{\ub} \in \reals^{|\PA|}$, $z_{k,l},z_{l,k} \in \qcone^3$ for $(k,l) \in \F$, $w_g \in \qcone^3$ for $g \in \G^{\text{quad}}$, $X \in \herm^n$, and $v \in \cmplx^n$. Notice that the constraints \eqref{eq:opf-reform-ls}-\eqref{eq:opf-reform-le} are all linear. We refer the reader to \citep{AHV:14} for a definition of the data matrices $Y_k, \widetilde Y_k, T_{k,l}$, and $\widetilde T_{k,l}$.

The only non-convex constraint in \eqref{eq:opf-reform} is the rank-1 condition \eqref{eq:opf-reform-rank1}. An SDR of \eqref{eq:opf-reform} is readily obtained by replacing \eqref{eq:opf-reform-rank1} by the positive semidefiniteness constraint $X \succeq 0$. The resulting SDR is a so-called cone linear program (CLP) that can be expressed as
\begin{align} \label{eq:cone-lp}
	\begin{array}{ll}
	\mbox{minimize}   & c^Tx \\
	\mbox{subject to} & Ax = b\\
                  	& x \in \K
	\end{array}
\end{align}
where $x$ is the vector of variables and the cone $\cone$ is a Cartesian product of three types of cones, \ie,
\[ \mathcal{K} = \reals_+^{n_l} \times \underbrace{\qcone^3 \times \cdots \times \qcone^3}_{n_q} \times \cone_{\mathrm{h}}^n.\]
Thus, the number of variables is $N = n_l + 3n_q +n^2$ where $n_l = 4|\G| + |\G^{\text{quad}}| + 2|\N| + 2|\PA|$ and $n_q = 2|\F|+|\G^{\text{quad}}|$, and the number of equality constraints is
$M = 4 |\N| + 2|\G| + 2 |\PA| + 3n_q$.

\subsection{Conversion}
The computational cost of solving \eqref{eq:cone-lp} with a general-purpose interior-point method becomes prohibitively large when $n$ is large: the cost of an interior-point iteration is at least $O(n^3)$.
Fortunately, the problem \eqref{eq:cone-lp} is generally very sparse in practice, and hence the conversion method of \citet{FKMN:01} may be used to rewrite \eqref{eq:cone-lp} as an equivalent CLP
\begin{align} \label{eq:cone-lp-convert}
	\begin{array}{ll}
	\mbox{minimize}   & \tilde c^T \tilde x \\
	\mbox{subject to} & \tilde A \tilde x = b\\
										& E \tilde x = 0 \\
                  	& \tilde x \in \widetilde \K
	\end{array}
\end{align}
with
\[ \widetilde \K = \reals_+^{n_l} \times \underbrace{\qcone^3 \times \cdots \times \qcone^3}_{n_q} \times \cone_{\mathrm{h}}^{r_1} \times \cdots \times \cone_{\mathrm{h}}^{r_m}. \]
The conversion essentially decomposes the cone $\cone_{\mathrm{h}}^n$ into a Cartesian product of a number of lower-dimensional cones $\cone_{\mathrm{h}}^{r_1}\times \cdots \times \cone_{\mathrm{h}}^{r_m}$ at the expense of a set of coupling constraints $E \tilde x = 0$. This reformulation of the problem can have a dramatic effect on the computational cost of solving the SDR of the ACOPF problem, and it effectively mitigates the $O(n^3)$ bottleneck that arises with the formulation \eqref{eq:cone-lp}. Moreover, the conversion technique often induces sparsity in the system of equations that define the search direction at each interior-point iteration, reducing the cost per iteration further if the solver can exploit this type of sparsity. The conversion technique was first applied to SDRs of the ACOPF problem by \citet{Jab:12}.

\subsection{Implementation} \label{implementation}
Before turning to our numerical experiments, we briefly outline our implementation \citep{opfsdr}. The code is written in Python and performs the following steps:
\begin{enumerate}
	\item Read case file and build the CLP \eqref{eq:cone-lp}.
	\item Apply conversion method: convert \eqref{eq:cone-lp} to \eqref{eq:cone-lp-convert}.
	\item Apply Hermitian-to-symmetric transformation: map $\cone_{\mathrm{h}}^{r_i}$ to $\cone_{\mathrm{s}}^{2r_i}$ for $i=1,\ldots,m$.
	\item Scale the problem data to improve conditioning.
\end{enumerate}
As part of the first step, we allow some preprocessing of the data: (i) slack variables $p_g$ (or $q_g$) for which $P_g^{\min} = P_g^{\max}$ (or $Q_g^{\min} = Q_g^{\max}$) may be eliminated, (ii) numerical proxies for infinity which are used to indicate the absence of limits (\eg, on generation) may be truncated, and (iii) a minimum resistance of transmission lines may be enforced. The Hermitian-to-symmetric transformation is a well known trick  that is only necessary because the solvers used in our experiments cannot directly handle cones of Hermitian positive semidefinite matrices; see \eg\ \citep{BoV:04}. The scaling of the problem data in step 4 yields an equivalent problem, and we found that for some solvers, this can reduce the computational time by roughly a factor of two; we briefly return to the topic of scaling in Section \ref{sec:discussion}.

\section{Results} \label{sec:results}
\subsection{Experiments}

To investigate the robustness and scalability of our methodology, we conducted a series of numerical experiments based on a collection of test cases from MATPOWER \citep{Zimmerman2011a} (which includes a number of test cases from \citep{josz2016}) and Power Grid Lib \citep{Pglib} with as many as $n= 70,\!000$ power buses; we have also included a synthetic case of the continental USA from the Electric Grid Test Case Repository \citep{Birchfield:2018} with $n= 82,\!000$ power buses. We excluded cases that are infeasible and cases with generator cost functions that are neither quadratic nor linear. For each test case, we set up a CLP formulation of the SDR and solved it using five different CLP solvers: MOSEK 8.1 \citep{mosek}, SeDuMi 1.3 \citep{Sturm1999}, SDPT3 4.0 \citep{sdpt3}, SCS 1.2.7 \citep{ocpb:16}, and CDCS 1.1 \citep{CDCS}. MOSEK, SeDuMi, and SDPT3 are interior-point methods whereas SCS and CDCS are first-order methods based on the alternating direction method of multipliers (ADMM).

To compare our methodology to an approach based on a modeling tool, we used SDPOPF \citep{Molzahn2013} from ``MATPOWER Extras'' to set up and solve an SDR of each case. SDPOPF uses YALMIP \citep{Lofberg2004} to set up the problem which is then solved numerically using one of several possible solvers: we used MOSEK in order to facilitate a fair comparison. Finally, to compare our approach to a nonlinear optimization approach, we used MATPOWER to set up and solve each case  with three different interior-point methods for nonlinear optimization: MIPS \citep{WMSZ+:07} from MATPOWER 6.1, IPOPT 3.12.9 \citep{Wachter2006} with PARDISO 6.0 \citep{Kourounis:2018}, and KNITRO 10.3.1 \citep{byrd2006}. These are all called via MATPOWER using its default initialization---the default is sometimes referred to as ``flat start'' since all voltages are set to 1 p.u. and the active power generation is set to the midpoint of its bounds. When successful, these solvers return a locally optimal solution that provides an upper bound on the optimal value in contrast to the SDR that provides a lower bound.

\subsection{Setup}
Using the implementation described in section \ref{implementation}, we processed the problem data before setting up the SDRs. Specifically, we truncated generator bounds larger than 50 times the base MVA. 
We remark that SDPOPF enforces a minimum transmission line resistance of $10^{-4}$; in the experiments, we do not enforce a minimum resistance in our SDR.

All experiments but those involving KNITRO were conducted on an HPC node with two Intel XeonE5-2650v4 processors (a total of 24 cores) and 240 GB memory. All experiments with KNITRO were conducted on different hardware (2.5 GHz Intel Core i5 CPU, 8 GB of memory) because of license restrictions. As a result, the KNITRO computation times that we report cannot be compared directly to those reported for the other solvers. All MATLAB-based solvers were used with MATLAB R2017b, and MOSEK was called through its Python interface in Python 3.6.3.  Finally, we modified the default solver options as follows: for SeDuMi, we raised the maximum number of iterations from 150 to 250; for SCS and CDCS, we limited the number of iterations to 20,000; for CDCS, we disabled ``chordalize'' and used the ``primal" solver since this allowed us to solve the most cases; for SCS, we used the direct solver; for SDPT3 we used a value of 400 for ``smallblockdim'' and changed the maximum number of iterations from 100 to 250.

\subsection{Robustness}
We start with an investigation of robustness. Table \ref{tab:status} contains a summary of return statuses for the different solvers for a total of 159 test cases. The column labeled ``success" refers to return values that indicate successful termination with an optimal or near optimal (global or local) solution. The ``failure'' column refers to return values that indicate some kind of error.
\begin{table}%[htb]
	\centering
	%\footnotesize
	\caption{Summary of return statuses by solver.}
	\label{tab:status}
	\begin{tabular}{l|rrr}
\toprule
Solver & Success & Max.\ iter. & Failure \\\midrule
MOSEK & 159 & 0 & 0\\ 
SeDuMi & 53 & 0 & 106\\SDPT3 & 52 & 0 & 107\\ 
SDPOPF & 128 & 0 & 31\\ 
CDCS & 146 & 13 & 0\\ 
SCS & 17 & 142 & 0\\ 
\midrule
IPOPT & 133 & 0 & 26\\ 
KNITRO & 145 & 0 & 14\\ 
MIPS & 116 & 0 & 43\\ 
\bottomrule
\end{tabular}
\end{table}
%We remark that MATPOWER fails to set up 9 cases from PGLIB\footnote{They have dispatchable loads where neither bound on the reactive power is zero ($Q_{\min}\ne0$ and $Q_{\max} \ne 0$).}, and hence never calls the solver in these cases. 
We remark that SDPOPF ignores phase angle constraints and fails in 31 cases because of a MATPOWER error; the solver is never called in these cases.

The results in Table \ref{tab:status} clearly demonstrate that the SDRs can be solved reliably using MOSEK: all cases were solved to optimality with MOSEK's default tolerances. In contrast, the nonlinear solvers IPOPT and KNITRO only succeed in roughly 85\% of the cases while MIPS succeeds in approximately 75\% of the cases. 
Note also that although CDCS solves all but one case, the accuracy and speed is poor compared to MOSEK as we show later in this section. 
Both SeDuMi and SCS succeed in less than 50\% of the cases.

%The fact that SDPOPF (which uses MOSEK as solver) only solves around 40\% of the cases emphasizes that the way the SDP problem is formulated and the solver is called affects the outcome significantly.
%It should be mentioned that many of the cases where SDPOPF failed is because angle difference constraints are not implemented.

\subsection{Accuracy}
\begin{figure*}%[htb]
	%\centering
	\subfloat[MOSEK]{
		\includegraphics[trim={{\myLeftTrim} {\myBotTrim} {\myRightTrim} {\myTopTrim}},clip=true,width=.48\linewidth]{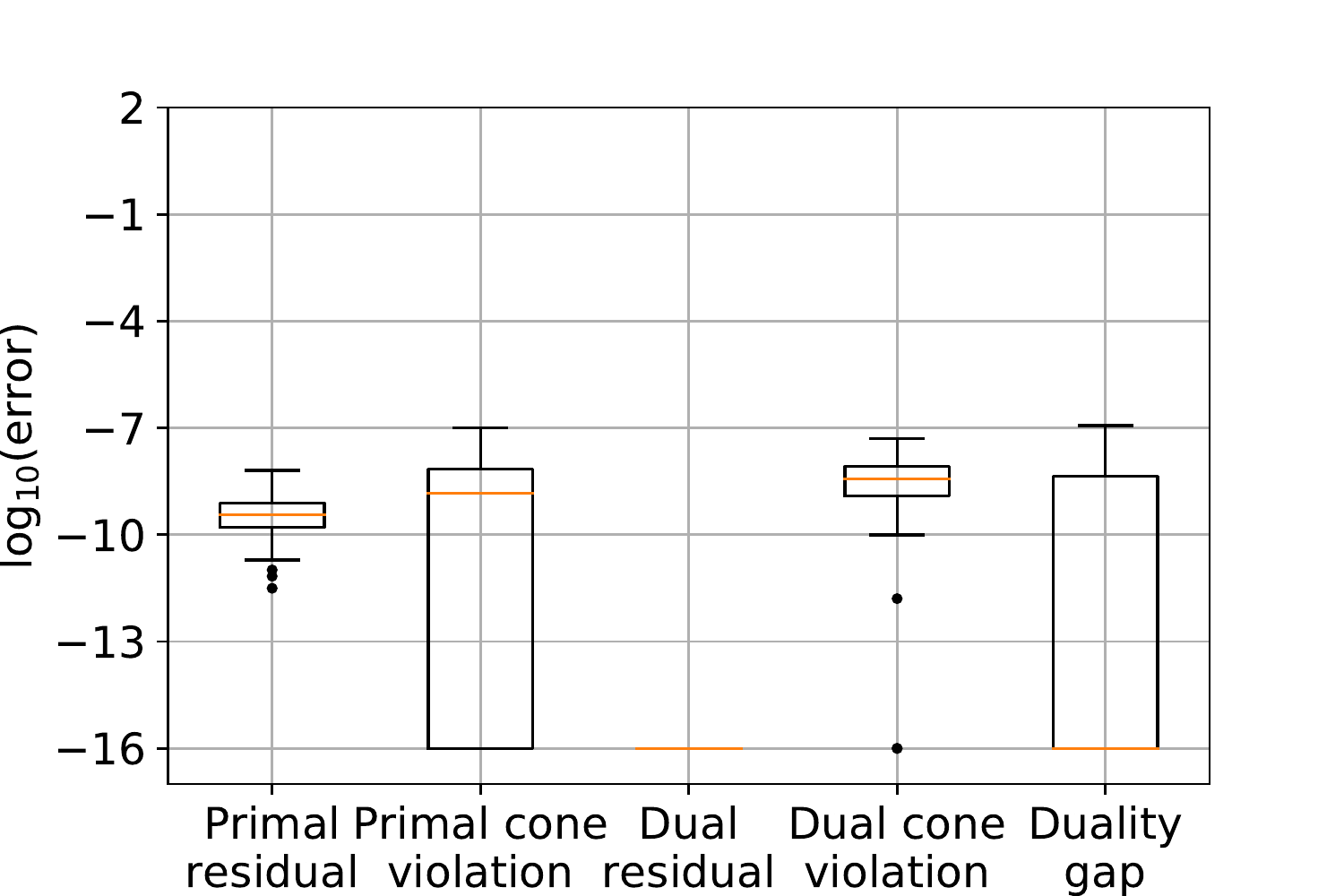}
		\label{fig:dimacsbox:mosek}
		}\hfill
	\subfloat[CDCS]{
		\includegraphics[trim={{\myLeftTrim} {\myBotTrim} {\myRightTrim} {\myTopTrim}},clip=true,width=.48\linewidth]{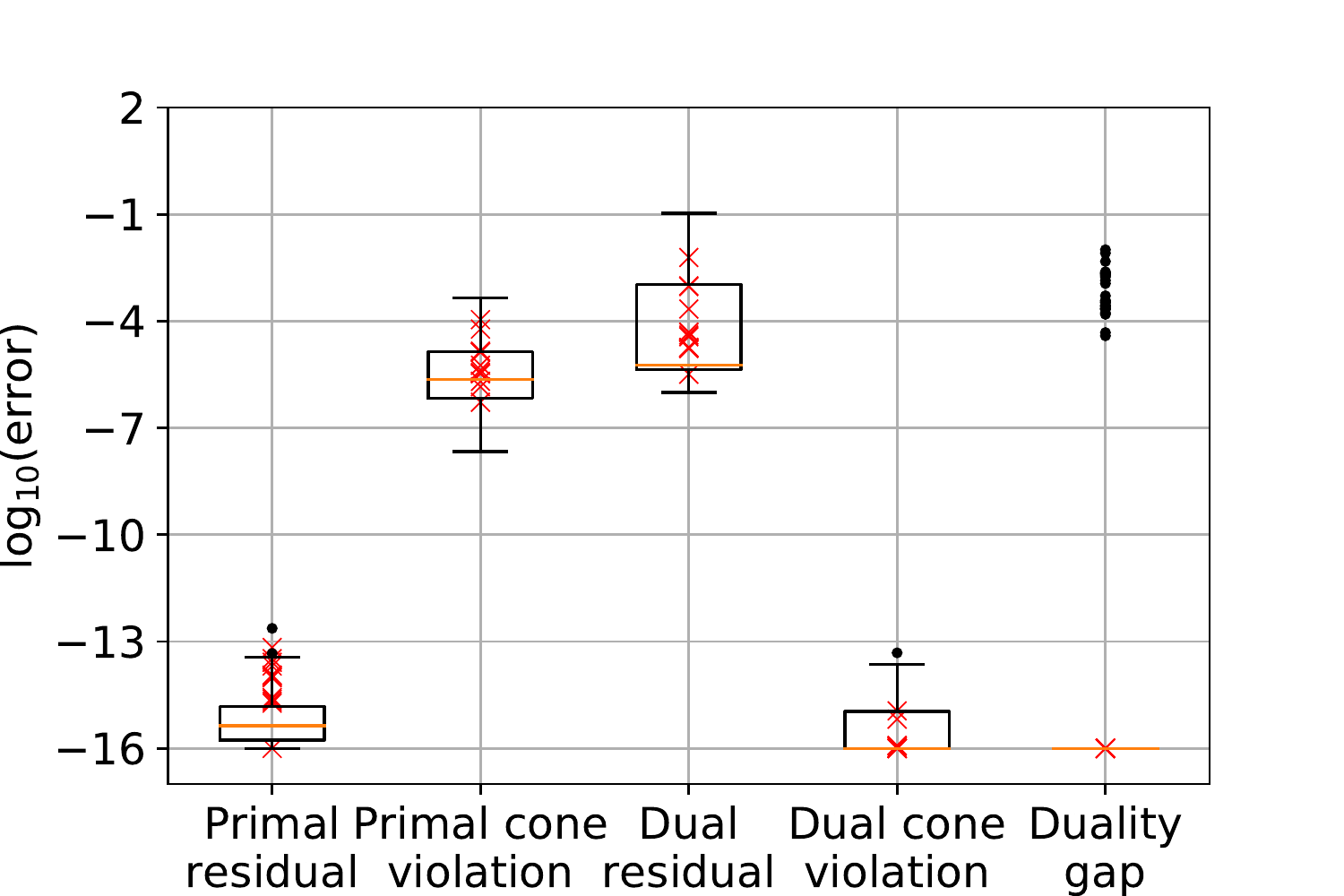}
		\label{fig:dimacsbox:cdcs}
		} \\[-1em]
	\subfloat[SeDuMi]{
		\includegraphics[trim={{\myLeftTrim} {\myBotTrim} {\myRightTrim} {\myTopTrim}},clip=true,width=.48\linewidth]{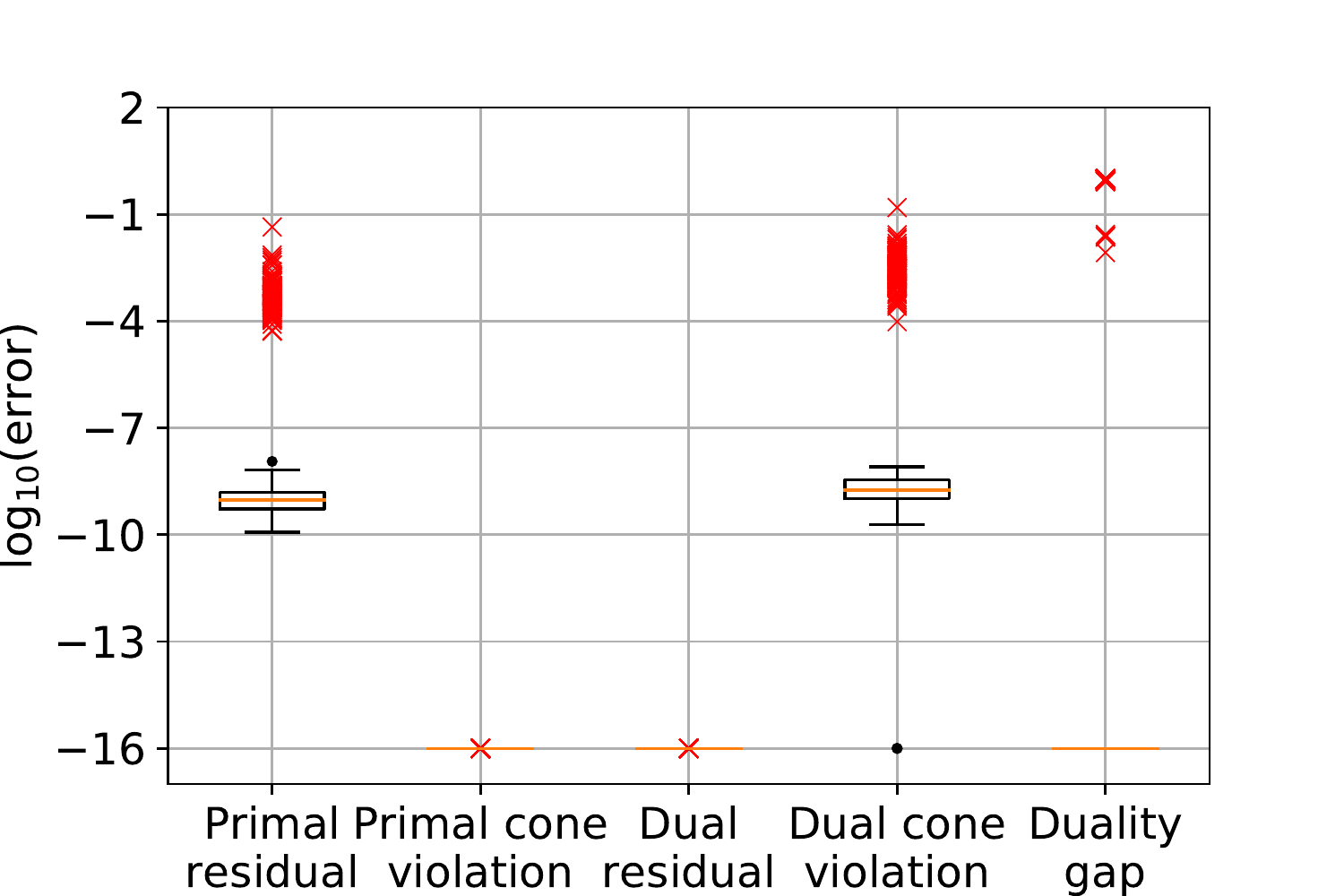}
		\label{fig:dimacsbox:sedumi}
		}\hfill
	\subfloat[SCS]{
		\includegraphics[trim={{\myLeftTrim} {\myBotTrim} {\myRightTrim} {\myTopTrim}},clip=true,width=.48\linewidth]{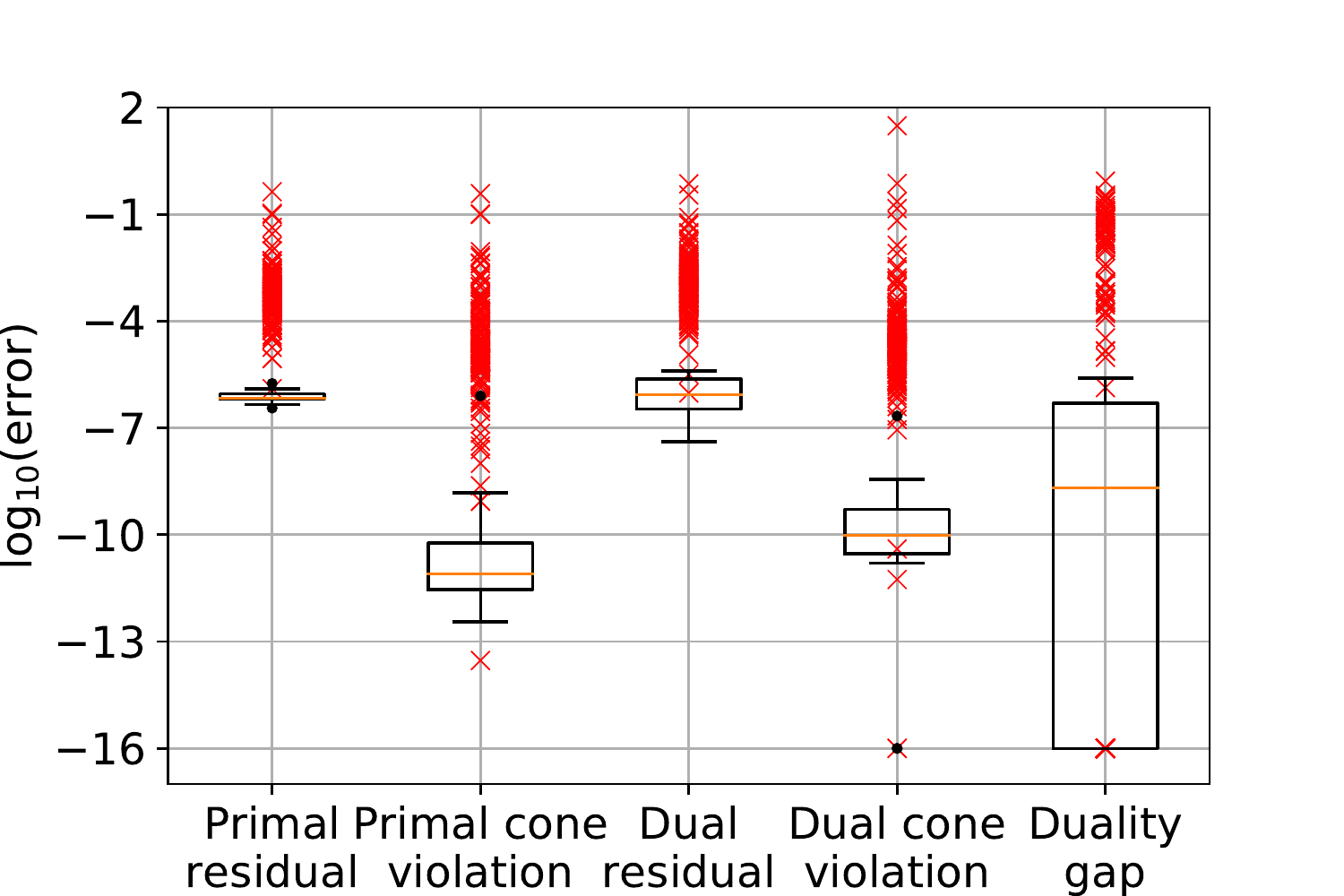}
		\label{fig:dimacsbox:scs}
		}\\[-1em]
	\subfloat[SDPT3]{
		\includegraphics[trim={{\myLeftTrim} {\myBotTrim} {\myRightTrim} {\myTopTrim}},clip=true,width=.48\linewidth]{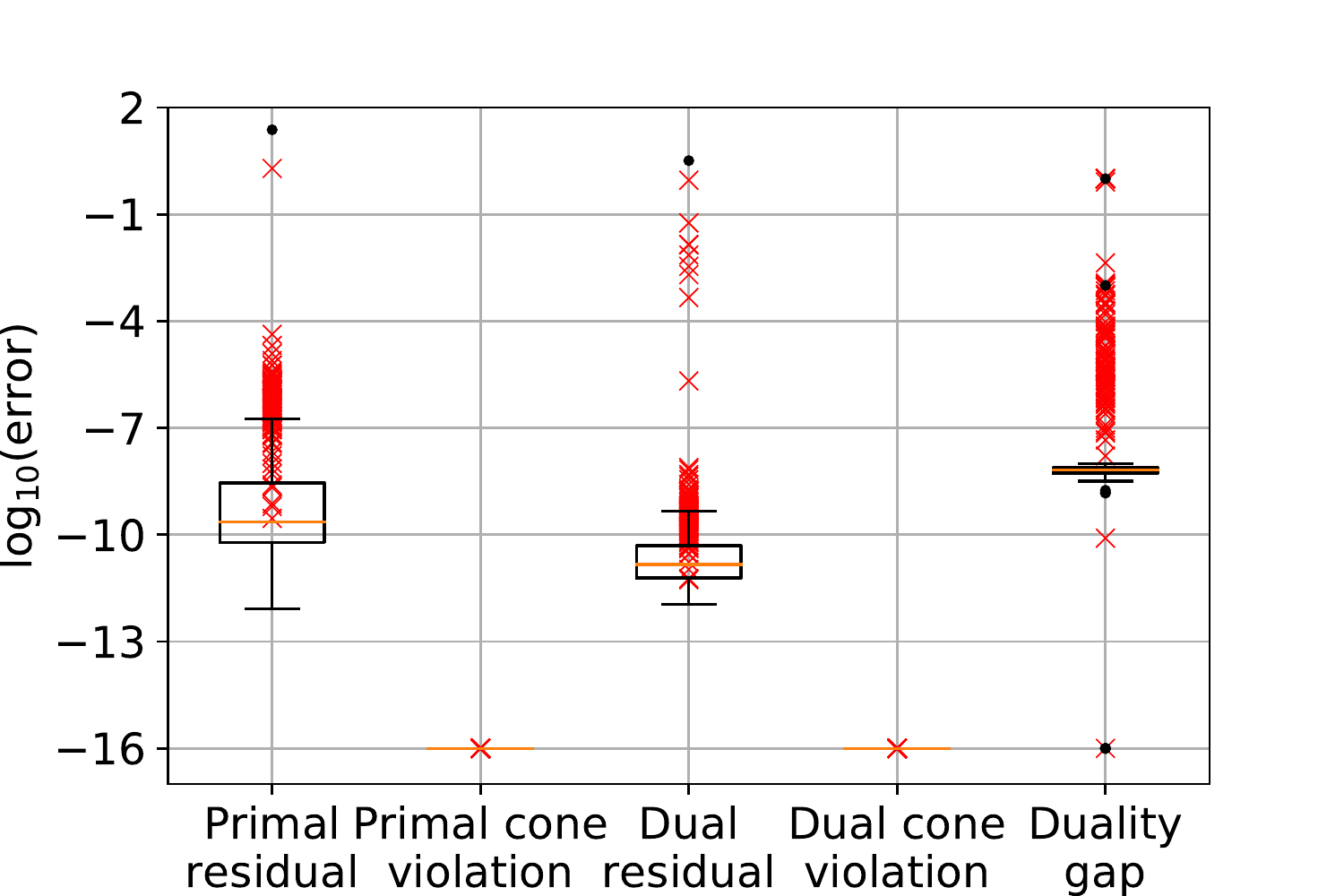}
		\label{fig:dimacsbox:sdpt3}
	}
	\caption{Box plots of logarithm of DIMACS errors. The red markers correspond to cases where the solver did not succeed. We note that in order to accommodate a logarithmic axis, we have replaced errors below $10^{-16}$ by this value.}
	\label{fig:dimacsbox}
\end{figure*}

We now compare the solutions returned by the five CLP solvers. Since the solvers have different tolerances (\ie, stopping criteria), we will compare the solvers based on the so-called ``DIMACS error measures" described in \citep{Mittelmann2003}. Roughly speaking, these are five relative error measures quantifying the primal residual norm, primal cone violation, dual residual norm, dual cone violation, and duality gap.
Fig.~\ref{fig:dimacsbox} summarizes the results in a box plot of the DIMACS measures for each solver (the smaller the error, the better).

MOSEK, shown in Fig.~\ref{fig:dimacsbox:mosek}, generally performs well with DIMACS errors below $10^{-7}$ in all cases. 
The SeDuMi errors, shown in Fig.~\ref{fig:dimacsbox:sedumi}, reveal that SeDuMi returns a high-accuracy solution whenever it succeeds; the same is true for SDPT3, shown in Fig.~\ref{fig:dimacsbox:sdpt3}. This suggests that the default tolerances may be too strict for all but the small cases. 
Both CDCS and SCS generally return solutions with larger errors, as shown in Fig.~\ref{fig:dimacsbox:cdcs} and \ref{fig:dimacsbox:scs}. This is to be expected since they are both first-order methods. While CDCS is relatively robust, it often terminates with sizable dual residuals which are indicative of low-accuracy solutions.

\subsection{Optimality Gap}
Next we investigate the objective values provided by the solvers. We limit our attention to MOSEK and the nonlinear solvers IPOPT, MIPS, and KNITRO. The nonlinear solvers provide an upper bound when they terminate at a feasible point. We define the best upper bound as
\begin{equation}
    \overline{f} = \min(f_{\text{IPOPT}},f_{\text{KNITRO}},f_{\text{MIPS}}),
\end{equation}
\ie, the minimum of the objective values provided by the three solvers (if a solver does not succeed, we define its objective value to be $\infty$).
Similarly, the SDR (MOSEK) provides a lower bound which we denote by $ \underline{f} = f_{\text{MOSEK}}$. The optimality gap may then be defined as
\begin{equation}
	\mathrm{gap} = \frac{\overline{f} - \underline{f}}{\overline{f}} \cdot 100\%.
\end{equation}
The gap is equal to $0$ if $ \underline{f} = \overline{f} $, implying that we have a globally optimal solution. On the other hand, if the gap is large, $ \overline{f} $ may be a poor local minimum and/or the SDR provides a weak lower bound $\underline{f}$.

We have made four tables listing objective values and optimality gap for all cases with more than 300 buses based on their origin: 
table \ref{tab:cases_matpower} contains cases from the MATPOWER library; 
table \ref{tab:cases_pglib} contains cases from PGLIB in typical operating conditions; 
table \ref{tab:cases_pglib_sad} contains cases from PGLIB with small phase angle difference constraints;
table \ref{tab:cases_pglib_api} contains cases from PGLIB with binding thermal limit constraints.
The cases are sorted by the number of buses in ascending order. Note that the optimality gap is undefined if none of the nonlinear solvers succeed. 
The optimality gap is close to zero in many cases and below 1\% in all but a handful of cases.

\subsection{Scalability}
\begin{figure}%[htb]
	\centering
	\includegraphics[width=0.9\columnwidth]{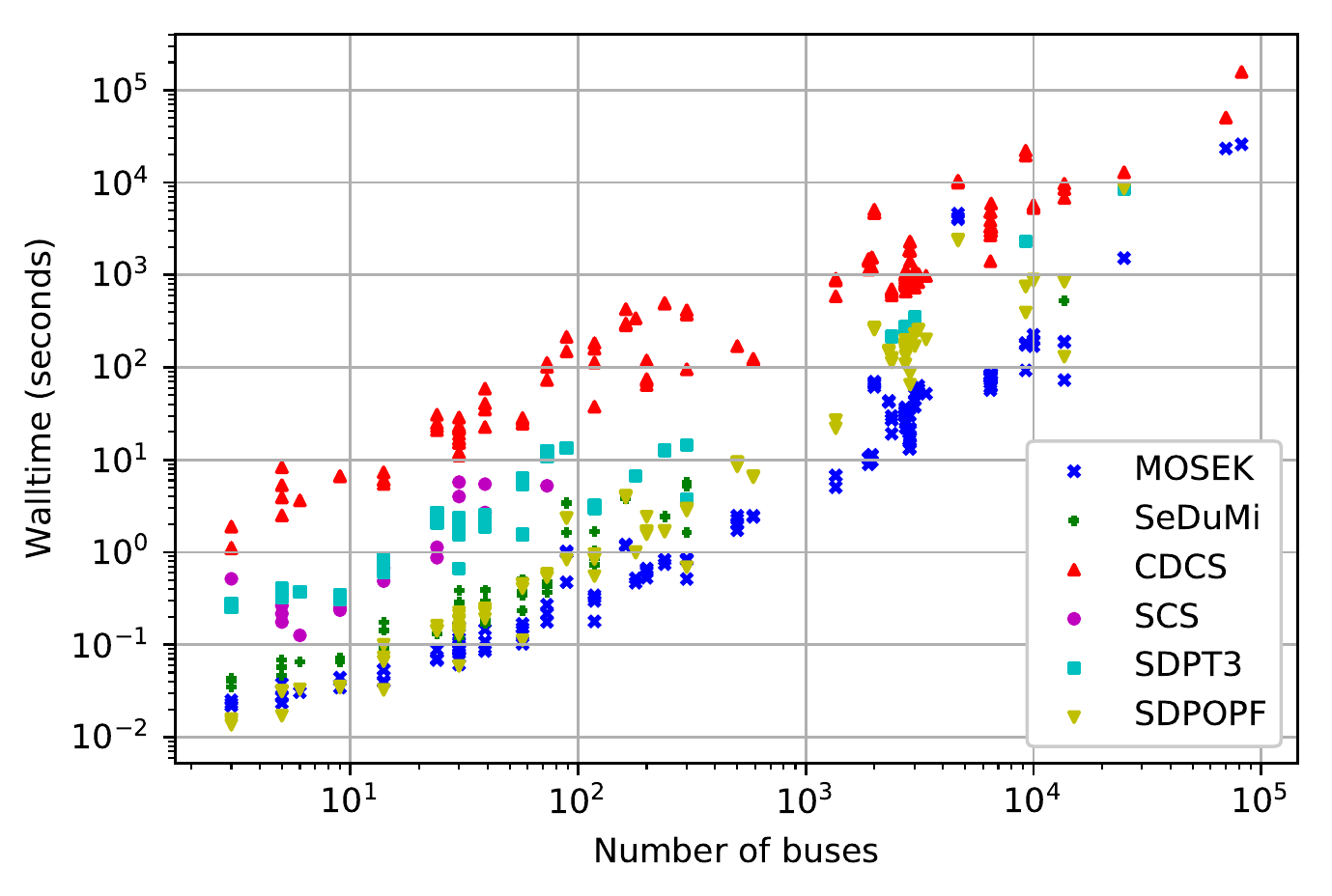}
	\caption{Scatter plot of the time used by the SDP solvers against the number of buses for successful cases. Note that SDPOPF solves an equivalent but different SDR.}
	\label{fig:timescat}
\end{figure}

We end this section by comparing the time required by each solver to solve the test cases. Fig.~\ref{fig:timescat} shows the time used by the SDP solvers compared to the number of buses in the case. To make a fair comparison, we report computation times without preprocessing, \ie, only the time required by the actual solver is recorded (we briefly discuss some considerations related to preprocessing in Section \ref{sec:discussion}).

MOSEK is generally the fastest. 
%Moreover, MOSEK is competitive with the nonlinear solvers (IPOPT, KNITRO, and MIPS) for small cases with a few hundred power buses. 
The difference between MOSEK and SDPOPF (which also uses MOSEK, but based on the problem formulation compiled by YALMIP) highlights that the formulation of the SDR may have a significant impact on the computation time as well as robustness. 
The striking difference between MOSEK and CDCS, both in terms of computation time and accuracy, makes it hard to justify the use of first-order methods for highly sparse problems like these.

In addition to cost function value and optimality gap, tables \ref{tab:cases_matpower}--\ref{tab:cases_pglib_api} list the computation times (excluding preprocessing) for MOSEK and the three solvers IPOPT, KNITRO, and MIPS. 
%The first block of cases are from MATPOWER, and the three subsequent blocks are cases from PGLIB with (i) typical operating conditions, (ii) binding thermal limit constraints (API), and (iii) binding phase angle difference constraints (SAD). 
MOSEK solves the SDR of all but one case with less than 25,000 buses in less than 10 minutes; the only exception is the case 4661\_sdet from PGLIB (all three operating conditions). Solving this problem takes MOSEK around 80 minutes. The longer computation time required to solve this case compared to other cases with a similar number of buses can in part be explain by looking at the chordal embedding of the network graph. The largest clique is of size 242 which is similar to the case with 82,000 buses (238) and around three times the size of all other cases with less than 25,000 buses. The case with 25,000 power buses is solved in approximately an hour by MOSEK, and the largest cases with 70,000 and 82,000 buses are solved in around seven hours. The nonlinear solvers are typically 5-20 times faster than MOSEK (they solve a different problem!), but they sometimes fail. The RTE cases from PGLIB appear to be particularly difficult for the nonlinear solvers: in some cases, none of the nonlinear solvers succeed, and the computation times are occasionally large compared to the general trend.

\begin{sidewaystable}[htb]
\centering
\caption{Cost, gap, and computation time for MATPOWER cases with more than 300 buses and without dispatchable loads. Failures are reported as `M' (max.\ iterations), `I' (termination at infeasible point), `N' (numerical error in solver). Times shown in \textcolor{red}{red} correspond to failures.}
\label{tab:cases_matpower}
\begin{tabular}{c|l|rrrr|r|dddd}
\toprule
&  & \multicolumn{4}{c|}{Cost} & \multicolumn{1}{c|}{Gap} & \multicolumn{4}{c}{Time (sec.)} \\
& Case & \multicolumn{1}{c}{IPOPT} & \multicolumn{1}{c}{KNITRO} & \multicolumn{1}{c}{MIPS} & \multicolumn{1}{c|}{MOSEK} &                     & \multicolumn{1}{c}{IPOPT} & \multicolumn{1}{c}{KNITRO} & \multicolumn{1}{c}{MIPS} & \multicolumn{1}{c}{MOSEK} \\
\midrule
\multirow{26}{*}{\rotatebox[origin=c]{90}{MATPOWER}}
 & ACTIVSg500 & 7.258e+04 & 7.258e+04 & 7.258e+04 & 7.105e+04 & 2.1\% & 1.6 & 0.7 & 0.6 & 2.0\\
 & 1354pegase & 7.407e+04 & 7.407e+04 & 7.407e+04 & 7.406e+04 & 0.0\% & 5.3 & 1.7 & 2.3 & 5.0\\
 & 1888rte & F & 5.980e+04 & F & 5.960e+04 & 0.3\% & \dcolcolor{red} 203 & 17.7 & \dcolcolor{red} 1.0 & 8.9\\
 & 1951rte & 8.174e+04 & 8.174e+04 & F & 8.173e+04 & 0.0\% & 19.9 & 9.7 & \dcolcolor{red} 1.2 & 9.7\\
 & ACTIVSg2000 & 1.228e+06 & 1.229e+06 & 1.228e+06 & 1.228e+06 & 0.0\% & 16.2 & 2.5 & 3.3 & 69.9\\
 & 2383wp & 1.869e+06 & 1.869e+06 & 1.869e+06 & 1.861e+06 & 0.4\% & 11.1 & 3.2 & 3.2 & 27.0\\
 & 2736sp & 1.308e+06 & 1.308e+06 & 1.308e+06 & 1.308e+06 & 0.0\% & 27.7 & 2.7 & 3.0 & 29.2\\
 & 2737sop & 7.776e+05 & 7.776e+05 & 7.776e+05 & 7.775e+05 & 0.0\% & 11.4 & 2.7 & 2.8 & 32.6\\
 & 2746wop & 1.208e+06 & 1.208e+06 & 1.208e+06 & 1.208e+06 & 0.0\% & 11.5 & 2.8 & 3.2 & 33.7\\
 & 2746wp & 1.632e+06 & 1.632e+06 & 1.632e+06 & 1.632e+06 & 0.0\% & 14.8 & 2.8 & 3.1 & 32.3\\
 & 2848rte & F & 5.302e+04 & F & 5.301e+04 & 0.0\% & \dcolcolor{red} 485 & 31.8 & \dcolcolor{red} 1.9 & 13.1\\
 & 2868rte & 7.979e+04 & 7.979e+04 & F & 7.979e+04 & 0.0\% & 80.7 & 10.7 & \dcolcolor{red} 2.4 & 13.3\\
 & 2869pegase & 1.340e+05 & 1.340e+05 & 1.340e+05 & 1.340e+05 & 0.0\% & 10.7 & 2.9 & 4.6 & 15.2\\
 & 3012wp & 2.592e+06 & 2.592e+06 & 2.592e+06 & 2.588e+06 & 0.1\% & 48.5 & 3.3 & 5.0 & 44.2\\
 & 3120sp & 2.143e+06 & 2.143e+06 & 2.143e+06 & 2.142e+06 & 0.0\% & 16.3 & 3.2 & 5.3 & 53.9\\
 & 3375wp & 7.412e+06 & 7.412e+06 & 7.412e+06 & 7.409e+06 & 0.0\% & 15.8 & 3.6 & 5.9 & 51.9\\
 & 6468rte & 8.683e+04 & 8.683e+04 & F & 8.682e+04 & 0.0\% & 56.3 & 23.5 & \dcolcolor{red} 10.7 & 56.7\\
 & 6470rte & 9.835e+04 & 9.835e+04 & F & 9.834e+04 & 0.0\% & 297 & 26.9 & \dcolcolor{red} 11.7 & 60.4\\
 & 6495rte & 1.063e+05 & 1.063e+05 & F & 1.061e+05 & 0.2\% & 198 & 25.0 & \dcolcolor{red} 8.7 & 61.7\\
 & 6515rte & 1.098e+05 & 1.098e+05 & F & 1.097e+05 & 0.1\% & 279 & 46.6 & \dcolcolor{red} 10.4 & 62.8\\
 & 9241pegase & 3.159e+05 & 3.159e+05 & 3.159e+05 & 3.158e+05 & 0.0\% & 208 & 131 & 19.8 & 92.8\\
 & ACTIVSg10k & 2.486e+06 & 2.486e+06 & F & 2.486e+06 & -0.0\% & 73.2 & 151 & \dcolcolor{red} 42.5 & 170\\
 & 13659pegase & F & 3.861e+05 & F & 3.861e+05 & 0.0\% & \dcolcolor{red} 1,\!221 & 72.3 & \dcolcolor{red} 3,\!953 & 72.9\\
 & ACTIVSg25k & 6.018e+06 & 6.018e+06 & F & 6.017e+06 & 0.0\% & 405 & 51.2 & \dcolcolor{red} 85.9 & 1,\!520\\
 & ACTIVSg70k & 1.644e+07 & 1.644e+07 & F & 1.644e+07 & 0.0\% & 896 & 199 & \dcolcolor{red} 140 & 23,\!343\\ \hline
 & SyntheticUSA & F & F & F & 2.017e+07 &  \textemdash  & \dcolcolor{red} 8,\!582 & \dcolcolor{red} 10,\!534 & \dcolcolor{red} 439 & 25,\!922\\
\bottomrule
\end{tabular}
\end{sidewaystable}
\begin{sidewaystable}[htb]
\centering
\caption{Cost, gap, and computation time for PGLIB cases in typical operating condition with more than 300 buses and without dispatchable loads. Failures are reported as `M' (max.\ iterations), `I' (termination at infeasible point), `N' (numerical error in solver). Times shown in \textcolor{red}{red} correspond to failures.}
\label{tab:cases_pglib}
\begin{tabular}{c|l|rrrr|r|dddd}
\toprule
&  & \multicolumn{4}{c|}{Cost} & \multicolumn{1}{c|}{Gap} & \multicolumn{4}{c}{Time (sec.)} \\
& Case & \multicolumn{1}{c}{IPOPT} & \multicolumn{1}{c}{KNITRO} & \multicolumn{1}{c}{MIPS} & \multicolumn{1}{c|}{MOSEK} &                     & \multicolumn{1}{c}{IPOPT} & \multicolumn{1}{c}{KNITRO} & \multicolumn{1}{c}{MIPS} & \multicolumn{1}{c}{MOSEK} \\
\midrule
\multirow{26}{*}{\rotatebox[origin=c]{90}{PGLIB}}
 & 500\_tamu & 7.258e+04 & 7.258e+04 & 7.258e+04 & 7.105e+04 & 2.1\% & 1.8 & 0.6 & 0.5 & 2.3\\
 & 588\_sdet & 3.816e+05 & 3.816e+05 & 3.816e+05 & 3.798e+05 & 0.4\% & 1.4 & 0.8 & 1.1 & 2.4\\
 & 1354\_pegase & 1.364e+06 & 1.364e+06 & 1.364e+06 & 1.356e+06 & 0.6\% & 5.0 & 2.1 & 2.7 & 6.7\\
 & 1888\_rte & 1.640e+06 & 1.565e+06 & F & 1.538e+06 & 1.7\% & 52.0 & 35.6 & \dcolcolor{red} 3.9 & 10.2\\
 & 1951\_rte & 2.375e+06 & F & F & 2.375e+06 & 0.0\% & 39.6 & \dcolcolor{red} 113 & \dcolcolor{red} 1.0 & 11.0\\
 & 2000\_tamu & 1.228e+06 & 1.228e+06 & 1.228e+06 & 1.228e+06 & 0.0\% & 17.0 & 2.6 & 3.7 & 65.2\\
 & 2316\_sdet & 2.257e+06 & 2.257e+06 & 2.257e+06 & 2.240e+06 & 0.7\% & 8.5 & 3.1 & 4.5 & 43.1\\
 & 2383wp\_k & 1.869e+06 & 1.869e+06 & 1.869e+06 & 1.861e+06 & 0.4\% & 11.8 & 4.0 & 3.2 & 27.4\\
 & 2736sp\_k & 1.308e+06 & 1.308e+06 & 1.308e+06 & 1.308e+06 & 0.0\% & 12.8 & 2.8 & 3.2 & 29.8\\
 & 2737sop\_k & 7.776e+05 & 7.776e+05 & 7.776e+05 & 7.775e+05 & 0.0\% & 10.9 & 3.0 & 3.0 & 31.1\\
 & 2746wp\_k & 1.632e+06 & 1.632e+06 & 1.632e+06 & 1.632e+06 & 0.0\% & 14.4 & 3.5 & 3.3 & 34.9\\
 & 2746wop\_k & 1.208e+06 & 1.208e+06 & 1.208e+06 & 1.208e+06 & 0.0\% & 12.5 & 3.6 & 3.6 & 36.7\\
 & 2848\_rte & 1.385e+06 & 1.385e+06 & F & 1.384e+06 & 0.0\% & 75.4 & 12.4 & \dcolcolor{red} 13.9 & 16.4\\
 & 2853\_sdet & F & 2.469e+06 & 2.469e+06 & 2.456e+06 & 0.5\% & \dcolcolor{red} 46.9 & 4.7 & 6.8 & 29.8\\
 & 2868\_rte & 2.260e+06 & 2.260e+06 & F & 2.260e+06 & -0.0\% & 43.1 & 18.7 & \dcolcolor{red} 18.0 & 17.4\\
 & 2869\_pegase & 2.605e+06 & 2.605e+06 & 2.605e+06 & 2.603e+06 & 0.1\% & 20.8 & 5.2 & 6.9 & 20.7\\
 & 3012wp\_k & 2.601e+06 & 2.601e+06 & 2.601e+06 & 2.597e+06 & 0.1\% & 21.0 & 4.0 & 5.4 & 50.7\\
 & 3120sp\_k & 2.146e+06 & 2.146e+06 & 2.146e+06 & 2.145e+06 & 0.0\% & 17.3 & 4.0 & 5.7 & 58.4\\
 & 4661\_sdet & F & 2.786e+06 & F & 2.768e+06 & 0.6\% & \dcolcolor{red} 554 & 8.9 & \dcolcolor{red} 19.7 & 4,\!184\\
 & 6468\_rte & F & 2.262e+06 & F & 2.252e+06 & 0.5\% & \dcolcolor{red} 1,\!476 & 76.5 & \dcolcolor{red} 3.5 & 82.5\\
 & 6470\_rte & F & F & F & 2.545e+06 &  \textemdash  & \dcolcolor{red} 1,\!434 & \dcolcolor{red} 31.5 & \dcolcolor{red} 7.6 & 77.7\\
 & 6495\_rte & 3.478e+06 & F & F & 2.966e+06 & 14.7\% & 489 & \dcolcolor{red} 153 & \dcolcolor{red} 44.3 & 77.1\\
 & 6515\_rte & F & 3.197e+06 & F & 2.992e+06 & 6.4\% & \dcolcolor{red} 1,\!403 & 61.4 & \dcolcolor{red} 39.3 & 84.0\\
 & 9241\_pegase & 6.775e+06 & 6.775e+06 & F & 6.770e+06 & 0.1\% & 399 & 27.9 & \dcolcolor{red} 75.4 & 175\\
 & 10000\_tamu & 2.486e+06 & 2.486e+06 & F & 2.486e+06 & -0.0\% & 98.7 & 113 & \dcolcolor{red} 40.1 & 195\\
 & 13659\_pegase & 1.078e+07 & 1.078e+07 & 1.078e+07 & 1.078e+07 & 0.0\% & 250 & 66.5 & 54.6 & 190\\
\bottomrule
\end{tabular}
\end{sidewaystable}
\begin{sidewaystable}[htb]
\centering
\caption{Cost, gap, and computation time for PGLIB cases with small angle differences with more than 300 buses and without dispatchable loads. Failures are reported as `M' (max.\ iterations), `I' (termination at infeasible point), `N' (numerical error in solver). Times shown in \textcolor{red}{red} correspond to failures.}
\label{tab:cases_pglib_sad}
\begin{tabular}{c|l|rrrr|r|dddd}
\toprule
&  & \multicolumn{4}{c|}{Cost} & \multicolumn{1}{c|}{Gap} & \multicolumn{4}{c}{Time (sec.)} \\
& Case & \multicolumn{1}{c}{IPOPT} & \multicolumn{1}{c}{KNITRO} & \multicolumn{1}{c}{MIPS} & \multicolumn{1}{c|}{MOSEK} &                     & \multicolumn{1}{c}{IPOPT} & \multicolumn{1}{c}{KNITRO} & \multicolumn{1}{c}{MIPS} & \multicolumn{1}{c}{MOSEK} \\
\midrule
\multirow{26}{*}{\rotatebox[origin=c]{90}{PGLIB SAD}}
 & 500\_tamu & 7.923e+04 & 7.923e+04 & 7.923e+04 & 7.322e+04 & 7.6\% & 2.7 & 0.7 & 0.6 & 2.5\\
 & 588\_sdet & 4.043e+05 & 4.043e+05 & 4.043e+05 & 3.814e+05 & 5.6\% & 2.1 & 0.8 & 1.2 & 2.5\\
 & 1354\_pegase & 1.365e+06 & 1.365e+06 & 1.365e+06 & 1.357e+06 & 0.6\% & 5.8 & 2.3 & 2.7 & 6.7\\
 & 1888\_rte & F & 1.640e+06 & F & 1.538e+06 & 6.2\% & \dcolcolor{red} 331 & 17.3 & \dcolcolor{red} 3.4 & 10.4\\
 & 1951\_rte & 2.383e+06 & F & F & 2.376e+06 & 0.3\% & 54.5 & \dcolcolor{red} 43.9 & \dcolcolor{red} 1.1 & 11.3\\
 & 2000\_tamu & 1.230e+06 & 1.230e+06 & 1.230e+06 & 1.229e+06 & 0.1\% & 31.2 & 3.2 & 3.9 & 61.4\\
 & 2316\_sdet & 2.257e+06 & 2.257e+06 & 2.257e+06 & 2.240e+06 & 0.7\% & 9.2 & 3.8 & 4.5 & 42.1\\
 & 2383wp\_k & 1.916e+06 & 1.916e+06 & 1.916e+06 & 1.905e+06 & 0.6\% & 15.8 & 3.9 & 3.4 & 29.8\\
 & 2736sp\_k & 1.329e+06 & 1.329e+06 & 1.329e+06 & 1.325e+06 & 0.4\% & 15.7 & 3.8 & 4.0 & 33.2\\
 & 2737sop\_k & 7.927e+05 & 7.927e+05 & 7.927e+05 & 7.859e+05 & 0.9\% & 15.8 & 4.2 & 3.8 & 33.4\\
 & 2746wp\_k & 1.667e+06 & 1.667e+06 & 1.667e+06 & 1.661e+06 & 0.4\% & 15.8 & 4.4 & 3.7 & 35.8\\
 & 2746wop\_k & 1.234e+06 & 1.234e+06 & 1.234e+06 & 1.226e+06 & 0.7\% & 16.7 & 3.9 & 3.8 & 36.9\\
 & 2848\_rte & F & F & F & 1.385e+06 &  \textemdash  & \dcolcolor{red} 558 & \dcolcolor{red} 106 & \dcolcolor{red} 13.0 & 16.4\\
 & 2853\_sdet & F & 2.495e+06 & 2.495e+06 & 2.458e+06 & 1.5\% & \dcolcolor{red} 87.5 & 5.4 & 6.5 & 30.7\\
 & 2868\_rte & F & F & F & 2.264e+06 &  \textemdash  & \dcolcolor{red} 41.3 & \dcolcolor{red} 21.7 & \dcolcolor{red} 2.7 & 17.4\\
 & 2869\_pegase & 2.620e+06 & 2.620e+06 & 2.620e+06 & 2.616e+06 & 0.2\% & 22.3 & 6.2 & 8.1 & 21.8\\
 & 3012wp\_k & 2.621e+06 & 2.621e+06 & 2.621e+06 & 2.610e+06 & 0.4\% & 21.9 & 5.2 & 6.0 & 52.7\\
 & 3120sp\_k & 2.176e+06 & 2.176e+06 & 2.176e+06 & 2.165e+06 & 0.5\% & 21.8 & 5.8 & 6.0 & 63.4\\
 & 4661\_sdet & F & 2.802e+06 & 2.802e+06 & 2.781e+06 & 0.7\% & \dcolcolor{red} 459 & 8.8 & 19.4 & 4,\!005\\
 & 6468\_rte & 2.262e+06 & 2.262e+06 & F & 2.252e+06 & 0.5\% & 434 & 45.2 & \dcolcolor{red} 3.5 & 80.9\\
 & 6470\_rte & F & F & F & 2.548e+06 &  \textemdash  & \dcolcolor{red} 2,\!637 & \dcolcolor{red} 1,\!519 & \dcolcolor{red} 16.1 & 78.6\\
 & 6495\_rte & F & 3.478e+06 & F & 2.966e+06 & 14.7\% & \dcolcolor{red} 2,\!458 & 54.2 & \dcolcolor{red} 22.7 & 82.0\\
 & 6515\_rte & F & F & F & 2.992e+06 &  \textemdash  & \dcolcolor{red} 1,\!234 & \dcolcolor{red} 143 & \dcolcolor{red} 6.1 & 79.3\\
 & 9241\_pegase & 6.920e+06 & 6.920e+06 & F & 6.823e+06 & 1.4\% & 142 & 29.1 & \dcolcolor{red} 89.5 & 183\\
 & 10000\_tamu & F & 2.486e+06 & F & 2.486e+06 & -0.0\% & \dcolcolor{red} 1,\!370 & 94.6 & \dcolcolor{red} 38.9 & 196\\
 & 13659\_pegase & 1.090e+07 & 1.090e+07 & 1.090e+07 & 1.082e+07 & 0.7\% & 187 & 42.7 & 55.5 & 188\\
\bottomrule
\end{tabular}
\end{sidewaystable}
\begin{sidewaystable}[htb]
\centering
\caption{Cost, gap, and computation time for heavily loaded PGLIB cases (i.e., binding thermal limits) with more than 300 buses and without dispatchable loads. Failures are reported as `M' (max.\ iterations), `I' (termination at infeasible point), `N' (numerical error in solver). Times shown in \textcolor{red}{red} correspond to failures.}
\label{tab:cases_pglib_api}
\begin{tabular}{c|l|rrrr|r|dddd}
\toprule
&  & \multicolumn{4}{c|}{Cost} & \multicolumn{1}{c|}{Gap} & \multicolumn{4}{c}{Time (sec.)} \\
& Case & \multicolumn{1}{c}{IPOPT} & \multicolumn{1}{c}{KNITRO} & \multicolumn{1}{c}{MIPS} & \multicolumn{1}{c|}{MOSEK} &                     & \multicolumn{1}{c}{IPOPT} & \multicolumn{1}{c}{KNITRO} & \multicolumn{1}{c}{MIPS} & \multicolumn{1}{c}{MOSEK} \\
\midrule
\multirow{18}{*}{\rotatebox[origin=c]{90}{PGLIB API}}
 & 500\_tamu & 4.034e+04 & 4.034e+04 & 4.034e+04 & 4.034e+04 & -0.0\% & 1.1 & 0.6 & 0.5 & 1.7\\
 & 588\_sdet & 4.996e+05 & 4.996e+05 & F & 4.983e+05 & 0.3\% & 1.4 & 0.8 & \dcolcolor{red} 0.4 & 2.4\\
 & 1888\_rte & F & 2.262e+06 & F & 2.259e+06 & 0.2\% & \dcolcolor{red} 672 & 8.0 & \dcolcolor{red} 4.1 & 10.6\\
 & 2000\_tamu & 1.288e+06 & 1.288e+06 & 1.288e+06 & 1.275e+06 & 1.0\% & 25.7 & 11.8 & 4.2 & 67.6\\
 & 2316\_sdet & 2.774e+06 & 2.774e+06 & 2.774e+06 & 2.758e+06 & 0.6\% & 9.1 & 3.8 & 4.1 & 43.2\\
 & 2383wp\_k & 2.791e+05 & 2.791e+05 & 2.791e+05 & 2.791e+05 & 0.0\% & 7.0 & 1.9 & 1.7 & 19.2\\
 & 2736sp\_k & 6.260e+05 & 6.260e+05 & 6.260e+05 & 6.097e+05 & 2.6\% & 14.7 & 4.3 & 3.5 & 30.7\\
 & 2737sop\_k & 3.587e+05 & 3.587e+05 & 3.587e+05 & 3.485e+05 & 2.8\% & 12.9 & 4.0 & 3.1 & 30.6\\
 & 2746wp\_k & 5.818e+05 & 5.818e+05 & 5.818e+05 & 5.818e+05 & 0.0\% & 7.6 & 2.6 & 2.2 & 22.1\\
 & 2746wop\_k & 5.117e+05 & 5.117e+05 & 5.117e+05 & 5.117e+05 & 0.0\% & 6.9 & 2.3 & 1.7 & 24.8\\
 & 3012wp\_k & 7.289e+05 & 7.289e+05 & 7.289e+05 & 7.289e+05 & 0.0\% & 10.8 & 2.6 & 5.2 & 37.5\\
 & 3120sp\_k & F & 9.696e+05 & 9.696e+05 & 8.818e+05 & 9.1\% & \dcolcolor{red} 43.8 & 5.3 & 5.8 & 57.5\\
 & 4661\_sdet & F & 3.343e+06 & F & 3.319e+06 & 0.7\% & \dcolcolor{red} 207 & 9.7 & \dcolcolor{red} 5.7 & 4,\!604\\
 & 6468\_rte & F & F & F & 2.718e+06 &  \textemdash  & \dcolcolor{red} 1,\!543 & \dcolcolor{red} 172 & \dcolcolor{red} 13.2 & 84.2\\
 & 6470\_rte & F & F & F & 3.174e+06 &  \textemdash  & \dcolcolor{red} 1,\!319 & \dcolcolor{red} 859 & \dcolcolor{red} 10.0 & 76.1\\
 & 6495\_rte & F & F & F & 3.735e+06 &  \textemdash  & \dcolcolor{red} 2,\!045 & \dcolcolor{red} 111 & \dcolcolor{red} 8.5 & 82.1\\
 & 6515\_rte & F & F & F & 3.657e+06 &  \textemdash  & \dcolcolor{red} 2,\!650 & \dcolcolor{red} 54.4 & \dcolcolor{red} 5.7 & 85.6\\
 & 10000\_tamu & 1.816e+06 & 1.816e+06 & F & 1.751e+06 & 3.6\% & 153 & 94.0 & \dcolcolor{red} 20.1 & 225\\
\bottomrule
\end{tabular}
\end{sidewaystable}

\iffalse
In Table \ref{tab:smallcases}, we see the results for the small-scale cases.
All cases are solved in less than 3 seconds by MOSEK and the optimality gap is generally small.
All cases except two (with binding thermal limits) are solved by IPOPT and the time used is generally less than MOSEK.

In Table \ref{tab:mediumcases}, we see the results for the medium-scale cases.
All cases are solved by MOSEK in less than 2 minutes.
IPOPT fails for cases from all libraries but it is fast for the cases that it does solve using less than 10 seconds for most cases; there are some (rte) cases where IPOPT uses more time than MOSEK.
The optimality gap is small (less than 1\%) for all cases except a few.

In Table \ref{tab:largecases}, we see the results for the large-scale cases.
Again MOSEK solves all cases and the time used is roughly 10 minutes for the largest case (scaling the cases the maximum time use by MOSEK was 3.5 minutes).
IPOPT solves around half the cases and the time used varies from very fast (less than 5 seconds) to immensely slow (more than 10 hours in the extreme).
The optimality gap is less than $ 1.5\% $ for all solved cases.
\fi

%{
%	\footnotesize
%	%\setlength\LTleft{-1in}
%	%\setlength\LTright{-1in}
%
%	\input{files/gapTable.txt}
%}

\section{Discussion} \label{sec:discussion}
The difference between our formulation of the SDR and the one constructed by SDPOPF via YALMIP shows that the problem formulation can have a significant impact on computation times and robustness. Our experiments demonstrate that an SDR of the ACOPF problem can be solved accurately and reliably with the right combination of problem formulation and solver. However, it is possible that the problem formulation can be further improved. For example, as mentioned in section \ref{implementation}, the conditioning of the problem may improve with some scaling of the constraints, and this, in turn, may reduce the number of iterations and/or the computation time. We have conducted some experiments in this direction, and our preliminary results show that using MOSEK, the solution time can roughly be cut in half; the geometric mean of the speed-up obtained by means of scaling was 1.9. Indeed, the solution time for the largest test case with 25,000 buses was reduced from about one hour to half an hour with MOSEK. We did not observe a similar improvement with scaling for the other solvers. Finally, we remark that scaling may affect stopping criteria, so care must be taken when comparing the accuracy of solutions obtained with and without scaling.

The QCR, proposed by \citet{coffrin2016}, provides a promising alternative to the SDR in that it is computationally cheaper and often as tight as the SDR (and in some cases even tighter). However, the findings reported in \citep{coffrin2016} only include SDRs of cases with less than 3,000 buses, and it is therefore unclear how the QCR and the SDR compare with respect to optimality gap for larger test cases. Moreover, the results pertaining to the SDR were obtained using an implementation based on SDPT3 and the modeling tool CVX, so the sizable gap between the two relaxations in terms of computational time will likely shrink if MOSEK and our problem formulation is used for the SDR.

%Recall that we used a minimum branch resistance of $10^{-5}$ in the experiments involving our problem formulation. We also tried solving the SDRs without enforcing a minimum branch resistance. The results obtained with MOSEK were very similar in both cases: MOSEK successfully solved all test cases with and without enforcing a minimum resistance.

The computation times reported in Section \ref{sec:results} did not include preprocessing time (\ie, the time required to construct the SDR). To give the reader an idea of the preprocessing workload, we remark that the construction of the SDR of the case with 25 thousand buses took approximately 25 seconds or approximately 1/60 of the time required to solve the SDR with MOSEK, and the geometric average of the ratio of the solution time to the preprocessing time for cases with more than 300 buses was approximately 13, \ie, preprocessing accounted for around 7\% of the total time on average. 
In contrast, YALMIP (via SDPOPF) required approximately 6 minutes to compile the case with 25,000 buses. 
Comparing the ratio of the preprocessing time for YALMIP to that of our approach, we found that the geometric average was approximately 13, \ie, on average it took 13 times longer with YALMIP. 
We note that our Python-based preprocessing code may be improved, \eg, by reimplementing critical parts of the code in C. 
In principle, the preprocessing time may be amortized if several problem instances with the same underlying power network need to be solved. 
However, this would require a symbolic chordal conversion of the problem such that the problem data can easily be updated or replaced. 

\section{Conclusion} \label{sec:conclusion}
SDR is a promising technique that may be used to compute useful global lower bounds on the optimal value of ACOPF problems. However, concerns about robustness and scalability have cast doubt on the practical usefulness of the technique. We have shown experimentally that the problem formulation can have a significant impact on both robustness and scalability. By constructing the SDR manually instead of using a modeling tool, we avoid problem transformations that incur significant overhead. Our numerical experiments establish that SDRs of a large collection of test cases can be solved reliably with MOSEK. Moreover, the time required to solve an SDR is typically within an order of magnitude of the time required by state-of-the-art nonlinear solvers such as KNITRO and IPOPT.

\bibliographystyle{apalike}
\bibliography{literature}

\end{document}